\documentclass[12pt]{amsart}
\usepackage{graphicx}
\usepackage{graphics}
\usepackage{amsmath}
\usepackage{amscd}
\usepackage{latexsym}
\usepackage{amssymb}
\usepackage{amsthm}
\usepackage{dsfont}
\usepackage{xcolor}

\usepackage{fancyhdr}
\usepackage[shortlabels]{enumitem}

\usepackage{tikz}
\usetikzlibrary{arrows,cd,matrix,positioning,decorations.pathmorphing}
\usepackage{setspace}
\tikzset{snake arrow/.style=
{->,
decorate,
decoration={snake,amplitude=.4mm,segment length=2mm,post length=1mm}},
}

\usepackage{blindtext}
\usepackage{scrextend}
\addtokomafont{labelinglabel}{\sffamily}

\usepackage{hyperref}  
\hypersetup{
  colorlinks   = true, 
  urlcolor     = blue, 
  linkcolor    = blue, 
  citecolor    = blue  
  }

\begin{document}

\textwidth 6.2in
\textheight 7.6in
\evensidemargin.75in
\oddsidemargin.75in

\newtheorem{Thm}{Theorem}
\newtheorem{Lem}[Thm]{Lemma}
\newtheorem{Cor}[Thm]{Corollary}
\newtheorem{Prop}[Thm]{Proposition}
\newtheorem{Rm}{Remark}
\newtheorem{Qu}{Question}
\newtheorem{Def}{Definition}
\newtheorem{Exm}{Example}

\def\a{{\mathbb a}}
\def\C{{\mathbb C}}
\def\A{{\mathbb A}}
\def\B{{\mathbb B}}
\def\D{{\mathbb D}}
\def\E{{\mathbb E}}
\def\R{{\mathbb R}}
\def\P{{\mathbb P}}
\def\S{{\mathbb S}}
\def\Z{{\mathbb Z}}
\def\O{{\mathbb O}}
\def\H{{\mathbb H}}
\def\V{{\mathbb V}}
\def\Q{{\mathbb Q}}
\def\Cn{${\mathcal C}_n$}
\def\CM{\mathcal M}
\def\CG{\mathcal G}
\def\CH{\mathcal H}
\def\CT{\mathcal T}
\def\CF{\mathcal F}
\def\CA{\mathcal A}
\def\CB{\mathcal B}
\def\CD{\mathcal D}
\def\CP{\mathcal P}
\def\CS{\mathcal S}
\def\CZ{\mathcal Z}
\def\CE{\mathcal E}
\def\CL{\mathcal L}
\def\CV{\mathcal V}
\def\CW{\mathcal W}
\def\CO{\mathcal O}
\def\CU{\mathcal U}
\def\O{{\mathbb O}}
\def\U{{\mathbb U}}

\def\IF{\mathbb F}
\def\IK{\mathcal K}
\def\IL{\mathcal L}
\def\IP{\bf P}
\def\IR{\mathbb R}
\def\IZ{\mathbb Z}

\title{On Shake Slice Knots}
\author{Selman Akbulut}
\author{Eylem ZEL\.{I}HA Yildiz}
\address{G\"{o}kova Geometry Topology Institute,  Mu\u{g}la, Turkiye}
\email{akbulut.selman@gmail.com}
\email{eylemzeliha@gmail.com}
\date{\today}
\begin{abstract} 
Here we discuss $r-$shake slice knots, and their relation to corks, we then prove that $0$-shake slice knots are slice.
\end{abstract}

\date{}
\maketitle

\setcounter{section}{-1}

\vspace{-.1in}

\section{Introduction and the main theorem} A knot $K\subset S^{3}$ is {\it slice} if it is the boundary of a properly imbedded smooth disk $D^{2}\subset B^{4}$. Let $K^{r}=B^{4}\smile h_{K}(r)$ be the $4$-manifold obtained by attaching  a $2$-handle to $B^{4}$ along the knot $K$ with framing $r$. Clearly $K$ is slice if and only if $K^{0}$ imbeds into $S^{4}$. To  see this, write $S^{4}$ as union of two $4$-balls $B^{4}_{-} \smile B^{4}_{+}$ glued along boundaries, then use the slice disk in $B^{4}_{+}$ to imbed $K^{0} = B_{-}^{4} \smile h_{K}(0)\subset S^{4}$ (compare \cite{km}).

\vspace{.1in}
                         
Clearly $K^{r}$ is  homotopy equivalent to $S^{2}$. We say that $K$ is {\it $r$-shake slice} if a generator of $H_{2}(K^{r})=\Z$ is represented by a smoothly imbedded $2$-sphere in $K^{r}$.
An {\it $r$-shaking of the knot $K$} is defined to be the link $\CL^{r} =\{K, K_{+},\ldots,K_{+}, K_{-},\ldots,K_{-}\}$, consisting of $K$ and an even number of oppositely oriented parallel copies of $K$ (pushed off by the framing $r$). Then by \cite{a2} and \cite{a1} we see that $K$  is r-shake slice if an $r$-shaking of $K$  bounds a disk $D_{n}$ with $2n$ holes for some $n\geq 0$ in $B^{4}$.  Now define $K^{r}_{n}= B^{4} \smile D_{n} \times B^{2}/\sim$, where we glue $D_{n} \times B^{2}$ to $B^{4}$ along $\partial D_{n} \times B^{2}$ to the framed link $\CL ^{r}\times D^{2}$ with $0$-framing. 

 \begin{figure}[ht] \begin{center}
 \includegraphics[width=.3\textwidth]{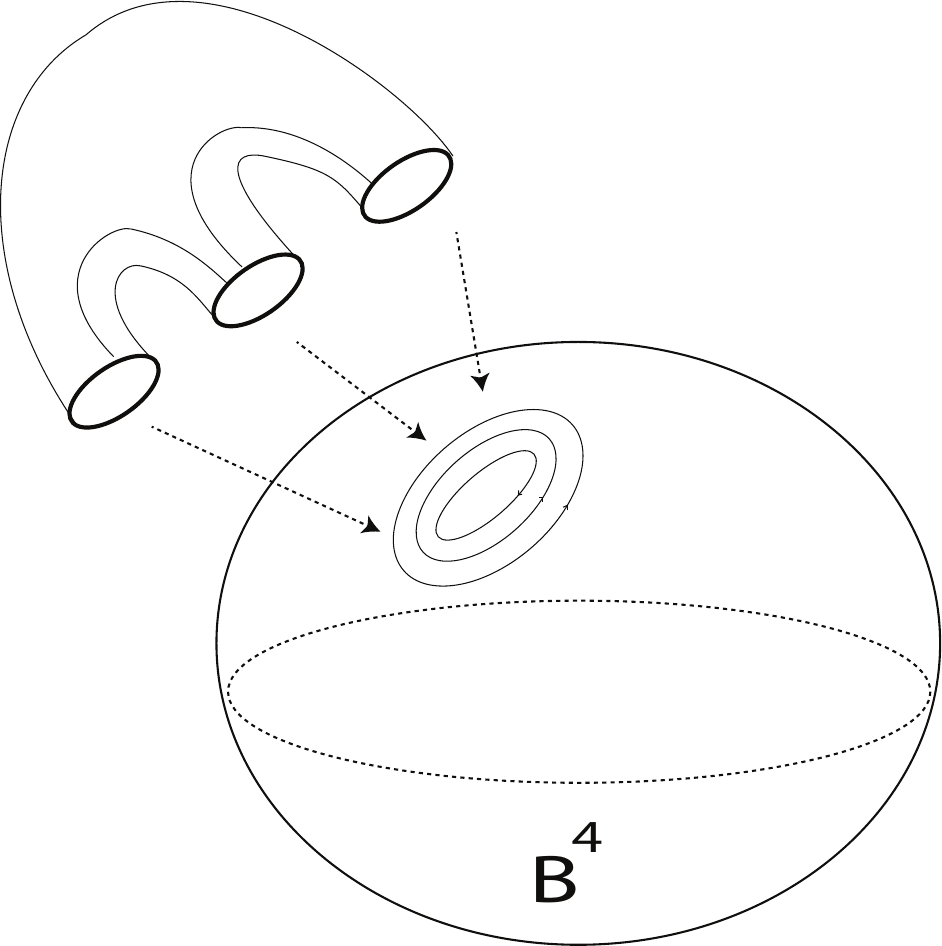} 
 \caption{$K_{n}^{0}$}
 \label{a1} 
\end{center}
 \end{figure}

\vspace{.1in}

It is known that, for each $r\neq 0$, there are $r$-shake slice knots which are not slice (Theorem 2.1 of \cite{a1}). Also not every knot is $r$-shake slice (e.g. Theorem 1 of \cite{a4}, and \cite{y}). By the same argument of the first paragraph above, we see that $K$ is $0$-shake slice if and only if $K_{n}^{0}$ imbeds into $S^{4}$ for some $n$. For example for $U_{n}^{0}\subset S^{4}$ for the unknot $U$. When $r\neq 0$ we also have: {\it $K$ is $r$-shake slice $\Rightarrow K^{r}_{n}\subset S^{4}$}. 

\vspace{.1in}

\begin{figure}[htbp] \begin{center}  
\includegraphics[width=.6 \textwidth]{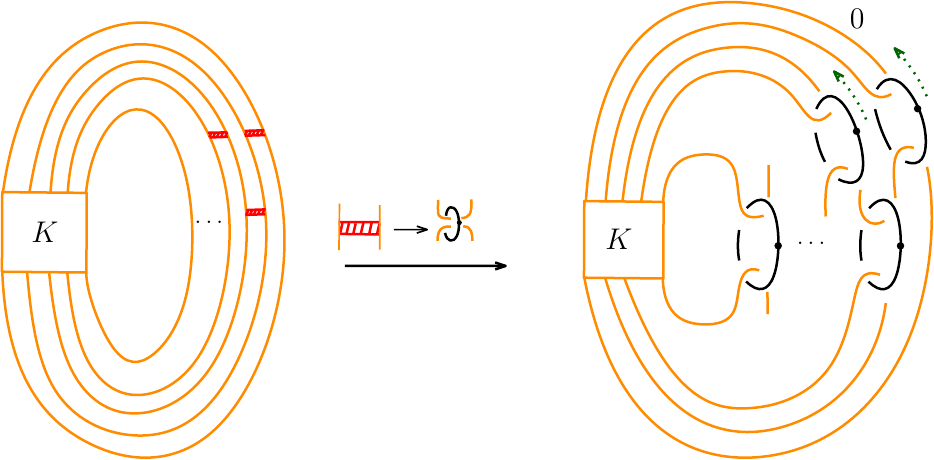}  
\caption{Handlebody of $K_{n}^{0}$}
\label{g1} 
\end{center}
\end{figure}

\begin{figure}[htbp]
 \begin{center}  
\includegraphics[width=.6 \textwidth]{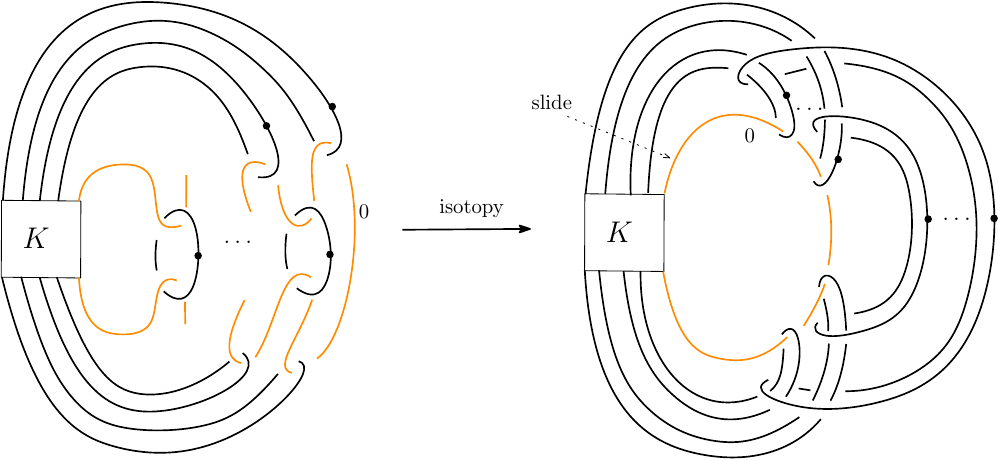}   
\caption{Handlebody of $K_{n}^{0}$}
\label{g2} 
\end{center}
\end{figure}

A handlebody of $K_{n}^{0}$ can be drawn by generalizing  the round handle attachment technique, applied  to $\CL^{0}$ (e.g. 3.2 of \cite{a1}),  i.e.  we attach $2$-handle with $0$-framing  to the knot obtained by band summing the link $\CL$ with $2n+1$ components $K\#(K_{+}\#...K_{+})\# (K_{-}\#...K_{-})$ so that each band going through a $1$-handle (circle with dot), as shown in Figure~\ref{g1}. After the obvious isotopies, from Figure~\ref{g1} we obtain Figure~\ref{g2}. 

\begin{figure}[htbp]  \begin{center}  
\includegraphics[width=.25 \textwidth]{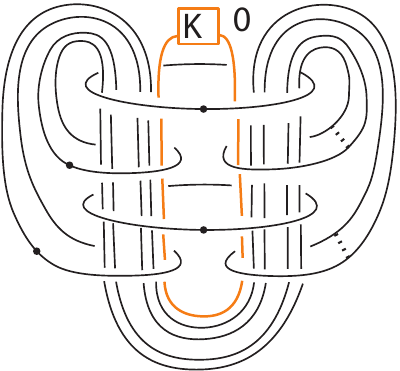}  
\caption{$K^{0}_{n}$,  (n=2)}
\label{g3}
\end{center}
\end{figure}

\begin{figure}[!h] 
	\centering 
	\begin{minipage}[t]{4.2cm} 
		\centering 
		\includegraphics[scale=0.5]{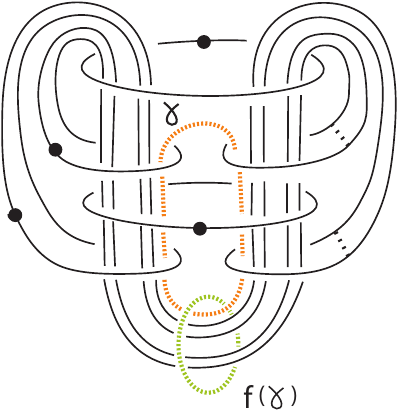} 
		\caption{$L_{n}^{0}$}
		\label{g3ab1}
	\end{minipage} 
	\hspace{.6cm} 
	\begin{minipage}[t]{5.2cm} 
		\centering 
		\includegraphics[scale=0.5]{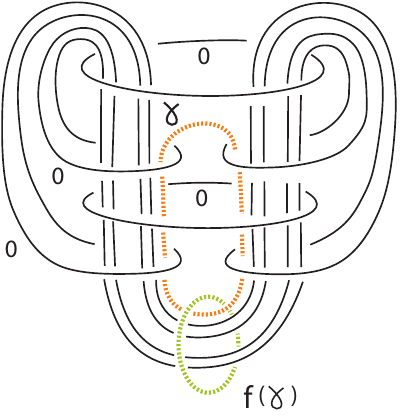} 
		\caption{$\small {S^{4}-L_{n}^{0}}$} 
		\label{g3ab2}
            \end{minipage} 
 \end{figure}

Notice, Figure~\ref{g2} shows $K^{0}_n$ is obtained by a single $2$-handle attachment to $\natural_{\,2n} S^1\times B^3$. Also note that $1$-handles of Figure~\ref{g2} are in the special position relative to the $2$-handle (as described in 1.4 of \cite{a1}), hence they can be slid over the $2$-handle of $K^{0}$ as indicated by Figure~\ref{g2} to obtain the handlebody  of $K_{n}^{0}$ in Figure~\ref{g3}. Which shows that $K_{n}^{0}$ is also obtained by a single $2$-handle attachment  $(\gamma \# K)^{0}$  to the manifold which we called $L_{n}^{0}$ in Figure~\ref{g3ab1}. 

 Figure~\ref{g3ab2} gives the complement of  the ``standard imbedding'' $L_{n}^{0}\subset S^{4}$, which is given by the obvious carved ribbons of Figure~\ref{g3ab1}. The choices of the carvings in Figure~\ref{g3ab1} determines the imbedding $L_{n}^{0}\subset S^{4}$, and $S^{4}-L_{n}^{0}$ in Figure~\ref{g3ab2}, but not vice versa. It is easy to see that, there is an involution $f:\partial L_{n}^{0}\to \partial L_{n}^{0}$ taking the loop $\gamma$ to $f(\gamma)$.
 
 \vspace{.1in}
 
\begin{figure}[htbp]  \begin{center}  
		\includegraphics[width=.38 \textwidth]{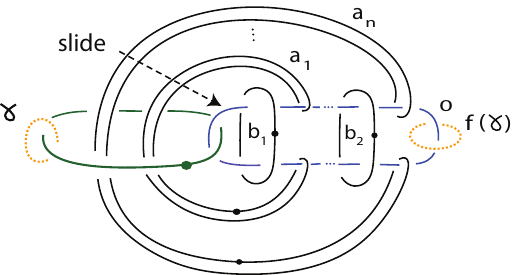}  
		\caption{$L_{n}^{0}$}
		\label{x}
	\end{center}
\end{figure}

 The following Proposition~\ref{L}  gives another handle description of   
 $L_{n}^{0}$ (Figure~\ref{x}), which consists of bunch of trivially carved discs from $B_{-}^{4}$(hence they are imbedded into $S^{4}$ uniquely), along with a $2$-handle. 
 
\newpage
 
  \begin{Prop} \label{L}
The handlebody of Figure~\ref{x} describes  $L_{n}^{0}$.  
  \end{Prop}
  
 \proof Slide the dotted circles $a_1, a_2,..,a_n$ over the $2$-handle (as described in Section 1.4 of \cite{a1}), then cancel the large $1$- and $2$- handle pair to get $L_{n}^{0}$  (Figure~\ref{g3ab1}). \qed
  
\vspace{-.1in}

\begin{Thm} \label{mainthm}
$0$-Shake slice knots are slice.
 \end{Thm}

\proof If a knot $K\subset S^{3}$ is $0$-shake slice, then $K_{n}^{0}\subset S^{4}$ for some $n$. Recall that$K_{n}^{0}$ is obtained by attaching a $2$-handle to $L_{n}^{0}$ along $\gamma \# K$ with $0$-faming, i.e. $K_{n}^{0}= L_{n}^{0}\smile \mbox{2-handle along}\; \gamma \# K \subset S^{4}$, 
which implies $\gamma \# K$ bounds a smooth disk $D$ in the complement $C: =S^{4}-L_{n}^{0}$.

\begin{figure}[htbp]  \begin{center}  
\includegraphics[width=.6\textwidth]{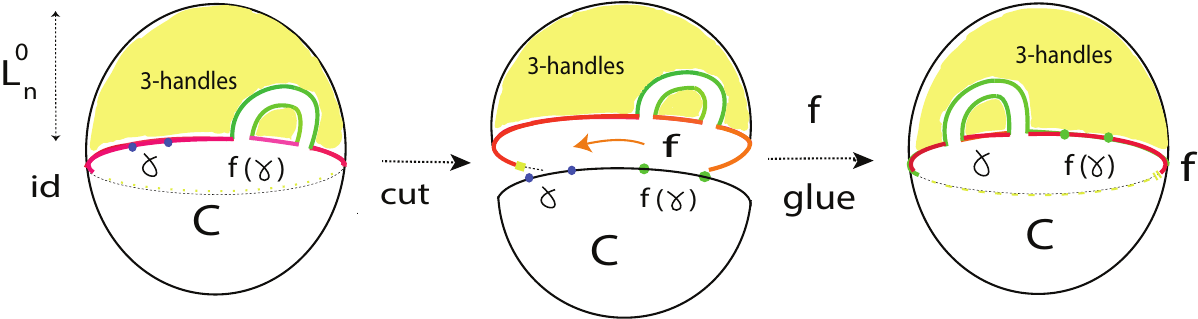}  
\caption{Cork twisting $S^{4}$. Carving a discs from $L_{n}^{0}$ is equivalent to attaching $2$-handle to its complement $C$}
\label{a6}
\end{center}
\end{figure}

From Figure~\ref{x} we see $f(\gamma)$ bounds a smooth disk in $L_{n}^{0}$ (the dual of the $2$-handle), which implies $f(\gamma\# K)$ bounds a singular disk $\bar{D}$ inside $L_{n}^{0}$ with $cone(K)$ singularity. If $\gamma \# K$ coincided with $f(\gamma \# K)$ we would be finished, because $D\smile_{\partial} \bar{D}$ would be a singular $S^{2}$  inside $S^{4}$ with $cone(K)$ singularity, which would imply $K^{0}\subset S^{4}$. Since this is not the case, we will complete the proof by performing a ``cork twisting operation'' to $S^{4}$ along $\partial L_{n}^{0}$, which throws $\gamma \# K$ to the $f(\gamma \# K)$ in $\partial L_{n}^{0}$, and then showing this doesn't change the smooth type of $S^{4}$. 

\vspace{.1in}

\begin{figure}[htbp]  \begin{center}  
		\includegraphics[width=.58 \textwidth]{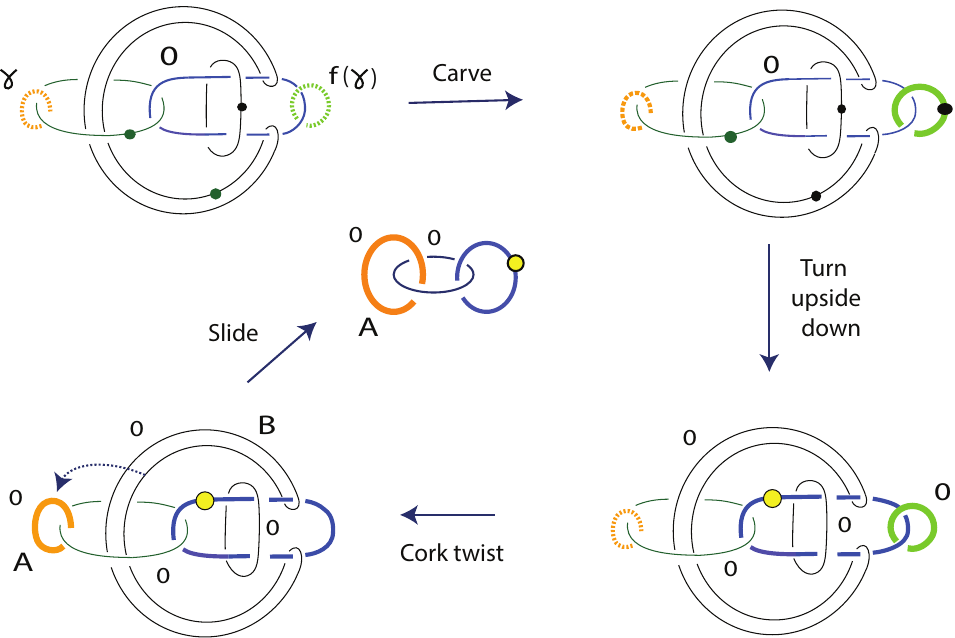}  
		\caption{Cork twisting $S^{4}$ along $L_{n}^{0}$}
		\label{y}
	\end{center}
\end{figure}

To see twisting $S^4$ along $\partial L_n^0$ via $f$ gives $S^4$ back, we recall 
$f(\gamma)$ bounds a smooth disk  $D\subset L_{n}^{0}$, and let $N(D)$ be its tubular neighborhood. Notice $L_{n}^{0}-N(D) \approx \natural_{2n+1} B^{3}\times S^{1}$, and upside down $\natural_{2n+1} B^{3}\times S^{1}$ are $3$-handles which are attached uniquely. Hence up to $3$-handles, the cork twisting $S^{4}$ along $L_n^0$ via $f$, is the operation described in  Figure~\ref{a6}
$$C\smile h^{2}_{f(\gamma)}\mapsto C\smile h^{2}_{\gamma}$$ 
To sum up, the first three steps of Figure~\ref{y} describe this operation:  $$L_{n}^{0}\longmapsto L_{n}^{0}-N(D) \longmapsto C\smile h^{2}_{f(\gamma)} \longmapsto C\smile h^{2}_{\gamma} $$

In Figure~\ref{y}  the circle with big yellow-dot denotes the complement of the $2$-handle of $L_{n}^{0}$. In the upside down picture, this corresponds to the complement of a slice $2$-disk  $\D$ (the core of the downside up $2$-handle).

\vspace{-.1in}

\vspace{-.1in} 

 Intersection patterns of this $2$-disk $\D$ with the $2$-handles of the third picture are dictated by the figure itself. Notice that this disk doesn't go over the $2$-handle $h^{2}_{\gamma}$ (=A) because $h^{2}_{\gamma}$ is attached later on. In the forth picture, $2$-handle denoted by B (and any piece of $\D$ which might be going over it) slides over the $2$-handle $h^{2}_{\gamma}$ to the other side of $h^{2}_{\gamma}$. Hence this carved disk $\D$ in the last picture doesn't go over the $2$-handle denoted by $A$; it lies entirely  in  $S{^2}\times B^{2}$ (after adding $3$-handles), and it has standard boundary $\partial B^{2}$. So by \cite{g} it is isotopic to the meridianal disk  p$\times B^{2}$. So it cancels the middle $2$-handle and results $B^{4}$ with a $2$-handle attached to it, with boundary $S^{2}\times S^{1}$. So by ``Property R" it is $S^{2}\times B^{2}$, which later gets cancelled by a $3$ handle to give $S^{4}$. \qed

\begin{figure}[htbp] 
  \begin{center}  
 \includegraphics[width=.6 \textwidth]{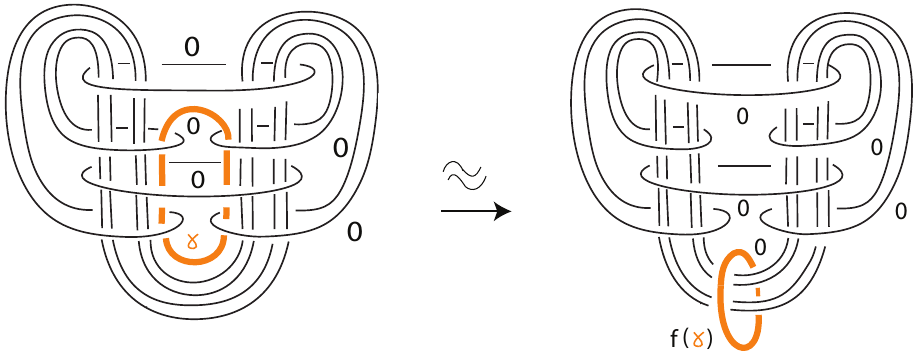}  
 \caption{Cork twisting $C\smile h^{2}_{\gamma} \mapsto C\smile h^{2}_{f(\gamma)}\approx S^{4}$}
 \label{b4d}
 \end{center}
 \end{figure}

\begin{Rm}
Start with the handlebody of $L_{n}^{0}$ in Figure~\ref{x}. Its $1$-handles $\#_{2n+1} B^{3}\times S^{1}$ can be imbeded into $S^{4}$ unique way, then by varying the imbedding of its $2$-handle in the complement we can otain many different imbeddings  $L_{n}^{0} \hookrightarrow S^{4}$.  By decomposing $S^{4}=B^{4}_{-}\smile B^{4}_{+}$, we can place $1$-handles as carved out $B^{4}_{-}$, and vary the imbedding of its $2$-handle in complement. For example, if we take the standard imbedding $L_{n}^{0}$ (Figure~\ref{g3ab1}) we get its complement $S^{4}-L_{n}^{0}$   (Figure~\ref{g3ab2}), and $f$ maps $h^{2}_{f(\gamma)} \mapsto  h^{2}_{\gamma}$. From Figure~\ref{b4d} we see that $C\smile h^{2}_{f(\gamma)} \approx C\smile h^{2}_{\gamma} \approx S^{4}$.
\end{Rm}

\section {Relation to corks}

The manifold $L_{n}^{0}$ has appeared in solution of many $4$-manifold problems. For example, it appeared in solution of Zeeman Conjecture \cite{a6}, it appeared in construction of diffeomorphic disk pairs $(B^{4}, D_{1})$, $(B^{4}, D_{2})$ that are not isotopic to each other rel boundary \cite{a7}. It  also appeard in $h$-cobordisms between homotopy $4$-spheres as  ``protocorks'' \cite{a8}. Here we will show that it can be used to generate new corks, from which we can construct new exotic $4$-manifolds (see \cite{a3}, \cite{a4}, \cite{a5}). 

 \begin{Thm} \label{cork}
  The diffeomorphism $f: \partial L^{0}_{n} \to \partial L^{0}_{n}$ can not extend to a self diffeomorphism of $ L^{0}_{n} \dashrightarrow  L^{0}_{n}$.
 \end{Thm}
 
 \proof We prove this by associating to $(L^{0}_{n}, f)$ a cork:  If f extended to a diffeomorphism $F: L^{0}_{n} \to  L^{0}_{n}$, then we could extend $F$ to a self diffeomorphism of the manifold $W_{n}$, obtained from $L^{0}_{n}$ by attaching $-1$ framed $2$-handles to $a_{1},.., a_{n}, b_{1},.., b_{n}$, since $f$ fixes these dotted circles.

\begin{figure}[htbp]  \begin{center}  
\includegraphics[width=.65\textwidth]{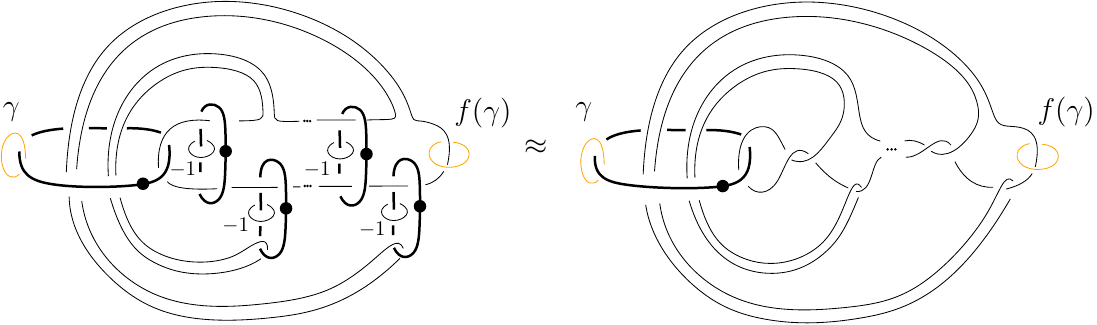}  
\caption{$W_{n}$}
\label{a}
\end{center}
\end{figure}

\begin{figure}[htbp]  \begin{center}  
\includegraphics[width=.34\textwidth]{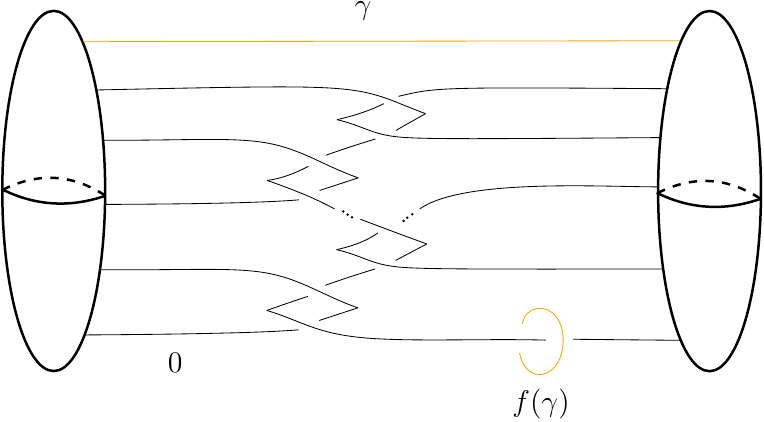}  
\caption{$W_{n}$}
\label{3}
\end{center}
\end{figure}

\noindent For simplicity, we continue to call  $F|_{\partial W_{n}} =f$.  Figure~\ref{a} is the handle picture of $W_{n}$, and Figure~\ref{3} is the Legendrian handle picture, from where we can check that it is a Stein manifold and calculate the TB invariant of $\gamma$, so by adjunction inequality it can not bound a smooth disk.  Since  $f(\gamma)$ bounds a smooth disk, the argument of Theorem 9.3 of \cite{a1} shows $(W_{n},f)$ is a cork and $f$ is a cork automorphism, so $F$ can't exist. $W_{n}$ is a variant of the ``positron cork'' of \cite{am}. \qed

 \begin{figure}[htbp]  \begin{center}  
\includegraphics[width=.17 \textwidth]{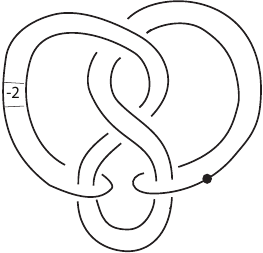}  
\caption{Exotic homotopy $S^{1}\times B^{3} $ rel boundary}
\label{g5a}
\end{center}
\end{figure}

 \newpage

 \begin{Rm}
$(L^{0}_{n}, f)$ is a basic universal object, which appears in construction of many exotic smooth $4$-manifolds.  Notice that the rel boundary exotic homotopy $S^{1}\times B^{3} $  of \cite{a6} (Figure~\ref{g5a}) can be obtained from  $L^{0}_{1}$ by attaching  $-1$ framed $2$-handle to $b_{1}$ in Figure~\ref{x},  which was carved out from the first cork introduced in \cite{a4}.   Notice the similarity between the construction of $K_{n}^{0}$  and construction of infinite order corks by altering concordances \cite{a7}.  
\end{Rm}

\begin{Rm}   However remote, there is some similarity between constructing handlebody of $K_{n}^{0}$ from $K$  and constructing ``spine manifolds''  from bounding manifolds, where shaking corresponds to taking disjoint union with spheres (Fact 3.2 of \cite{ak}).
 \end{Rm}



\begin{thebibliography}{99999}



\bibitem[A1]{a1} S. Akbulut, {\em $4$-Manifolds}, OxfordUP (2016) 

\bibitem[A2]{a2} S. Akbulut, {\em On 2-dimensional homology classes of 4-manifolds }, Math. Proc. Camb. Phil. Soc. 82 (1977), 99-106.

\bibitem[A3]{a3} S. Akbulut, {\em Corks }, https://arxiv.org/pdf/2406.15369


\bibitem[A4]{a4} S. Akbulut,  {\em A Fake compact contractible 4-manifold}, Journ. of Diff. Geom. 33, (1991), 335-356.

\bibitem[A5]{a5} S. Akbulut, {\em  An exotic 4-manifold}, J. Differ. Geom. 33 (1991) 357–361.

\bibitem[A6]{a6} S. Akbulut, {\em A solution to a conjecture of Zeeman}, Topology, vol.30, no.3, (1991), 513-515.


\bibitem[A7]{a7} S. Akbulut, {\em Corks and exotic ribbons in $B^4$} European Journal of Mathematics (2022)


\bibitem[A8]{a8} S. Akbulut, {\em On 4-dimensional smooth Poincare Conjecture}.\\
https://arxiv.org/pdf/2209.09968


\bibitem[AK]{ak} S. Akbulut and H.King  {\em The topology of real algebraic sets with isolated singularities}, Ann. of Math. 113 (1981), 425-446.

\bibitem[AM]{am} S. Akbulut and R. Matveyev  {\em Convex decomposition of 4-manifolds} IMRN International Math Research Notices, 1998, No.7 371-381.

\bibitem[G]{g} D. Gabai {\em The $4$-dimensional light bulb theorem} \\ https://arxiv.org/pdf/1705.09989.pdf

\bibitem[KM]{km} R. Kirby and P. Melvin, {\em Slice Knots and Property R},  Invent. Math. 45, 57- 59 (1978). 

\bibitem[Y]{y} K. Yasui, {\em Corks, exotic 4-manifolds and knot concordance}. \\arXiv:1505.02551v2.

\end{thebibliography}
\end{document}